\documentclass[12pt]{amsart}
\usepackage{}
\usepackage{mathrsfs}
\usepackage{dsfont,bbm,color}
\usepackage{amssymb,amsmath,amsthm}
\numberwithin{equation}{section}
\newtheorem{theorem}{Theorem}[section]
\newtheorem{lemma}[theorem]{Lemma}
\newtheorem{remark}[theorem]{Remark}
\newtheorem{example}[theorem]{Example}
\newtheorem{definition}[theorem]{Definition}

\begin{document}

\newcommand{\bbmu}{{\mbox{\boldmath$\mu$}}}
\newcommand{\gS}{\Sigma}
\newcommand{\bbgS}{{\mbox{\boldmath$\gS$}}}
\newcommand{\diag}{{\mbox{diag}}}
\newcommand{\gl}{{\lambda}}
\newcommand{\bbV}{{\bf V}}

\newcommand{\us}{\underline{s}}
\newcommand{\um}{\underline{m}}
\newcommand{\E}{\mathbb{E}}
\newcommand{\var}{\mathbb{V}ar}
\newcommand{\cov}{\mathbb{C}ov}
\renewcommand{\P}{\mathbb{P}}
\renewcommand{\a}{\alpha}
\renewcommand{\b}{\beta}
\newcommand{\Be}{\pmb\beta_0}
\renewcommand{\k}{\kappa}
\renewcommand{\t}{\tau}
\newcommand{\ep}{\varepsilon}
\newcommand{\tep}{\tilde\varepsilon}
\newcommand{\hep}{\hat\varepsilon}
\newcommand{\noin}{\noindent}
\newcommand{\De}{\Delta}
\newcommand{\la}{\lambda}
\newcommand{\La}{\Lambda}
\newcommand{\de}{\delta}
\renewcommand{\th}{\theta}
\newcommand{\si}{\sigma}
\newcommand{\Ga}{\Gamma}
\newcommand{\ga}{\gamma}
\newcommand{\pmu}{\pmb\mu}
\newcommand{\Si}{\pmb\Sigma}
\renewcommand{\o}{\omega}
\renewcommand{\O}{\Omega}
\newcommand{\ol}{\overline}
\newcommand{\bbI}{{\mathbf I}}
\newcommand{\bbC}{{\mathbb C}}
\newcommand{\bbR}{{\mathbb R}}
\newcommand{\cC}{{\mathcal C}}
\newcommand{\bbN}{{\mathbb N}}
\newcommand{\cD}{{\mathcal D}}
\newcommand{\cG}{{\mathcal G}}
\newcommand{\cI}{{\mathcal I}}
\newcommand{\cJ}{{\mathcal J}}
\newcommand{\cK}{{\mathcal K}}
\newcommand{\cQ}{{\mathcal Q}}

\newcommand{\frB}{\mathfrak{B}}
\newcommand{\frm}{\mathfrak{m}}
\newcommand{\frt}{\mathfrak{t}}
\newcommand{\cL}{{\mathcal L}}
\newcommand{\cH}{{\mathcal H}}
\newcommand{\cR}{{\mathcal R}}
\newcommand{\cT}{{\mathcal T}}
\newcommand{\cU}{{\mathcal U}}
\newcommand{\cV}{{\mathcal V}}
\newcommand{\vp}{\varpi}
\newcommand{\trvp}{\varpi^{tr}}
\newcommand{\evp}{\varpi^{\E}}

\newcommand{\A}{{\bf A}}
\newcommand{\C}{{\bf C}}
\newcommand{\B}{{\bf B}}
\newcommand{\D}{{\bf D}}
\newcommand{\e}{{\bf e}}
\newcommand{\F}{{\bf F}}
\newcommand{\G}{{\bf G}}
\renewcommand{\H}{{\bf H}}
\newcommand{\I}{{\bf I}}
\newcommand{\J}{{\bf J}}
\newcommand{\K}{{\bf K}}
\renewcommand{\L}{{\bf L}}
\newcommand{\M}{{\bf M}}
\newcommand{\N}{{\bf N}}
\newcommand{\Q}{{\bf Q}}
\renewcommand{\O}{{\bf O}}
\renewcommand{\r}{{\bf r}}
\renewcommand{\S}{{\bf S}}
\newcommand{\T}{{\bf T}}
\newcommand{\U}{{\bf U}}
\newcommand{\V}{{\bf V}}
\newcommand{\W}{{\bf W}}
\newcommand{\X}{{\bf X}}
\newcommand{\x}{{\bf x}}
\newcommand{\Y}{{\bf Y}}
\newcommand{\Z}{{\bf Z}}
\newcommand{\bz}{{\bf z}}
\newcommand{\bv}{{\bf v}}
\renewcommand{\(}{\left(}
\renewcommand{\)}{\right)}
\newcommand{\lj}{\left|}
\newcommand{\rj}{\right|}
\newcommand{\lb}{\label}
\newcommand{\no}{\nonumber\\}

\newcommand{\hi}{H\"{o}lder's inequality}
\newcommand{\eqD}{\stackrel{\mathcal{D}}{=}}
\newcommand{\toD}{\stackrel{\mathcal{D}}{\to}}

\title[]{STRONG REPRESENTATION OF WEAK CONVERGENCE\\[0.5em] \small (D\lowercase{edicate to the memory of} W\lowercase{enqi} L\lowercase{iang and} W\lowercase{im}  V\lowercase{ervaat}) }
\author{Zhidong Bai}
\address{KLASMOE and School of Mathematics \& Statistics, Northeast Normal University, Changchun, 130024, P.R.C.}
\email{baizd@nenu.edu.cn}

\author{Jiang Hu}
\address{KLASMOE and School of Mathematics \& Statistics, Northeast Normal University, Changchun, 130024, P.R.C..}
\email{huj156@nenu.edu.cn}
\thanks{Z. D. Bai was partially supported by CNSF 11171057. J. Hu was partially supported by NSFC 11201175.}


\maketitle

\begin{abstract}
Skorokhod's representation theorem  states that if on a  Polish  space, there is defined a weakly convergent sequence of probability measures $\mu_n\stackrel{w}\to\mu_0,$ as $n\to \infty$, then there exist a probability space $(\Omega, \mathscr F, P)$ and a sequence of random elements $X_n$ such that $X_n\to X$ almost surely and $X_n$ has the distribution function $\mu_n$, $n=0,1,2,\cdots$.
In this paper, we shall extend the Skorokhod representation theorem to the case where if there are a sequence of separable metric spaces $S_n$, a sequence of  probability measures $\mu_n$ and a sequence of measurable mappings $\varphi_n$ such that $\mu_n\varphi_n^{-1}\stackrel {w}\to\mu_0$, then there exist a probability space $(\Omega,\mathscr F,P)$ and $S_n$-valued random elements $X_n$ defined on $\Omega$, with distribution $\mu_n$ and such that
$\varphi_n(X_n)\to X_0$ almost surely. In addition, we  present several applications of our result including some results in random matrix theory, while the original Skorokhod representation theorem is not applicable.
\end{abstract}

\section{ Introduction and main result.}

Skorohod in 1956 \cite{Skorokhod56L} established his famous strong representation theorem that if $P_n$, $n = 1, 2,\dots$ and $P_0$ are  probability measures on a complete and separable
metric topological space $S$ (known as Polish space) such that $P_n\stackrel{w}\to P_0$, as $n\to\infty$, then there exist a probability space $(\Omega, \mathscr F, {P})$ and a sequence of measurable random elements $X_n$, $n=0,1,\dots,$ such that $P_n$ is the distribution of $X_n$ (i.e. $X_n\sim P_n$) and that $X_n\to X_0$ almost surely (a.s). Later Dudley in \cite{Dudley68D}
successively removed the completeness condition. Skorokhod representation theorem deeply reveals the relationship between convergence in distribution and strong convergence. Consequently, this  theorem  has been a strong mathematical tool of intense research of weak convergence for more than six decades. It realizes a convergent sequence of probability measures as distributions of a convergent random elements defined on the metric space $\Omega$ and thus it serves as a theoretical basis for the continuous-mapping approach in stochastic process. Skorokhod's representation theorem has many extensions and applications, more details can be found in \cite{BlackwellD83E,BertiP10S,Sethuraman02S,Whitt02S}.

In probability theory and mathematical physics, a random matrix is a matrix-valued random variable.
In late 1980's,  the spectral theory, which studies the distribution of the eigenvalues as the size of the matrix goes to infinity,
plays a very important role in the research of large dimensional random matrices. Apparently, if  these random variables are defined in different metric spaces and the dimensions of these matrices go to infinity,  we couldn't use the original Skorohod's  representation theorem directly, it had been found that Skorokhod theorem is not convenient to use when dealing such problems and thus Bai and Liang in \cite{Bailiang} extended Skorokhod theorem to a sequence of probability measures $\mu_n$ defined on a sequence of Polish spaces $S_n$ such that $\mu_n\varphi_n^{-1}\stackrel{w}\to\mu_0$, where $\varphi_n$ is a sequence of measurable mappings. Later, Wim Vervaat (when he was an associate editor of Annals of Probability) simplified the proof of the theorem and extended it to the non-completeness cases (see \cite{Bailv}). However due to the pass-away of Professor Vervaat in 1994, the paper has not been formally published yet. As a memorial to Professors Vervaat and Liang (1930-2007) and due to its great applicability in random matrix theory, we would like to formally publish the result in the present paper. In addition, to illustrate its powerful applications, we shall present several examples in which the original Skorokhod theorem is not applicable.

Throughout this paper,
equality in distribution is denoted by $\stackrel{\mathcal{D}}{=}$, convergence in distribution by $\stackrel{\mathcal{D}}{\to}$,  convergence in probability by $\stackrel{{p}}{\to}$, and weak convergence of probability measure by $\stackrel{w}\to$. Our main result of this paper is as follows.

\begin{theorem}\lb{th1}
Let $S_n$ for $n =0, 1,2, \dots $ be a sequence of separable metric spaces and let $\varphi_n$ for $n =1,2, \dots$  be a sequence of
measurable mappings from $S_n$ into $S_0$. Suppose that $\mu_n$ is a probability measure defined on the Borel field $\mathscr B_{S_n}$
generated by open subsets of $S_n$ for $n =0, 1,2, \dots$ and that
$\mu_n\varphi_n^{-1} \stackrel{w}{\to}\mu_0$, then there exist $S_n$-valued random elements $X_n$ for $n =0, 1,2, \dots$ defined on a common
probability space $(\Omega, \mathscr F, P)$ and such that $X_n\sim \mu_n$ in $(S_n,\varrho_n)$ for $n =0, 1,2, \dots$  and $\varphi_n(X_n) \to  X_0$ with probability one (w.p.1) as $n\to\infty$. Here $\varrho_n$ is the metric on $S_n$.
\end{theorem}
\begin{remark}
When $S_n \equiv S_0$ (separable) and $\varphi_n$ are identities for all $n\ge 1$, then the above theorem specializes to Dudley's
 variant of Skorohod's representation theorem.
\end{remark}

The rest of this paper is organized as follows.  In Section 2 we give the proof
of Theorem 1.1 and some following results.  Some technique lemmas are given in
Section 3 and some applications of Theorem 1.1 are given in Section 4. In Section 4
we present some results of random matrix theory as well.

\section{Proof of Theorem \ref{th1}}
In this section we give the proof of Theorem \ref{th1} and some following results.
Our proof of the present theorem amounts to the construction of a special metric space $T$
to which Dudley's theorem (Theorem 1 in \cite{Dudley68D}) can be applied. Notice that all  statements  involving  $n$  are  supposed  to  hold  for  $n  = 0, 1,2, \dots, $  unless  restricted
explicitly;  limit  statements  without  explicit  tendency  hold  as  $n  \to\infty$.

\textbf{Proof of Theorem \ref{th1}.}  We first consider the case where all $\varphi_n$ are continuous mappings. Let  $T$  be  the
disjoint  union  of  all $S_n$. Define space-indicator $s$ from  $T$  onto $\{0, 1 , 2, \dots, \}$   by $s(x)  =  n$  if  $x  \in  S_n$.   Set $ \varphi_{0}:=id_{S_0}$
and  define  $\varphi  : T  \to S_0$
by  $\varphi (x)  : = \varphi_{s(x)}(x)$.   Let  $\varrho_n$  denote  the  metric  on $S_n$.   Let $\ep_n$  be  positive  for  $n  >0$,  decreasing  to  0  as
$n  \to\infty$,  and  set  $\ep_0:=0$.    We  now define a  metric on  $T$ by
\begin{align}\lb{eq1}
    \de(x,y):=\varrho_0(\varphi(x),\varphi(y))+\left\{
                                                \begin{array}{ll}
                                                  \ep_{s(x)}\wedge\varrho_{s(x)}(x,y), & \hbox{if $s(x)=s(y)$;} \\
                                                  \ep_{s(x)}\vee\ep_{s(y)}, & \hbox{if $s(x)\neq s(y)$.}
                                                \end{array}
                                              \right.
\end{align}

Our  first task is to verify  that  $\de$  is  indeed  a  metric.   Obviously,  $\de(x ,y)  =  \de  (y,x)\ge 0$  and $\de(x ,x)  =  0$.    If  $\de  (x ,y) =
0$,  then  $s (x)  =  s (y)$  and  $\varrho_{s(x)}(x ,y)  =   0$,  thus  $x  = y$.  The  triangle  inequality can  be  verified  separately  for
both terms  on the  right-hand side  of (\ref{eq1}). In fact, we denote the two terms on the right hand side of (\ref{eq1}) by $\de_1(x,y)$ and
$\de_2(x,y)$, respectively. Then, it is obvious that for any $x,y,z\in T$,
$$
\de_1(x,z)=\varrho_0(\varphi(x),\varphi(z))\le \varrho_0(\varphi(x),\varphi(y))+\varrho_0(\varphi(y),\varphi(z))
\le \de_1(x,y)+\de_1(y,z)
$$
because $\varrho_0$ is a metric in $S_0$.

Now, if $ s(x)=s(y)=s(z)$, then
$$
\begin{cases}
\de_2(x,y)+\de_2(y,z)
=\varrho_{s(x)}(x,y)+\varrho_{s(x)}(y,z)\ge
 \varrho_{s(x)}(x,z)\ge \de_2(x,z)\cr
 \ \ \ \ \ \ \ \ \ \ \ \mbox{ if } \ep_{s(x)}\ge \max(\varrho_{s(x)}(x,y),\varrho_{s(x)}(y,z))\cr
\de_2(x,y)+\de_2(y,z)
\ge \ep_{s(x)}\ge \de_2(x,z)
\mbox{ if } \ep_{s(x)}< \max(\varrho_{s(x)}(x,y),\varrho_{s(x)}(y,z))\cr
\end{cases}
$$
If $ s(x)=s(y)\ne s(z)$, then we have
$$
\de_2(x,y)+\de_2(y,z)\ge \de_2(y,z)=\de_2(x,z).
$$
Symmetrically, if $ s(x)\ne s(y)= s(z)$,  we have
$$
\de_2(x,y)+\de_2(y,z)\ge \de_2(x,y)=\de_2(x,z).
$$
The last case, if $ s(x)\ne s(y) \ne s(z)$, then we have
$$
\de_2(x,y)+\de_2(y,z)=\ep_{s(x)}\vee \ep_{s(y)}+ \ep_{s(y)}\vee \ep_{s(z)} \ge \ep_{s(x)}\vee \ep_{s(z)}
\ge \de_2(x,z).
$$
Note that the last inequality above is an equality if $s(x)\ne s(z)$. Consequently, we have proved that the triangular inequality holds for the function $\de=\de_1+\de_2$. Thus $\de$ is a metric.

Our next task is to verify that the metric space $(T,\de)$ is separable. Let $Q_n$ be the countable dense subset of $S_n$ with respect to the topology generated by $\varrho_n$, $n=0,1,2,\cdots$.
Then by Lemma \ref{lem1}, for any open set $B$ of the metric space $(T,\de)$, $B\cap S_n$ is an open subset of $(S_n,\varrho_n)$ and hence contains an element of $Q_n$. Therefore, any open
 subset of $(T,\de)$ contains an  element of $Q=\cup_{n=0}^\infty Q_n$. Therefore, the metric $(T,\de)$ is separable.

Finally, to apply Dudley's theorem, we need to define the probability measures $\tilde\mu_n$ by $\tilde\mu_n(B)=\mu_n(B\cap S_n)$ for all Borel sets $B$ of $(T,\de)$. In this definition, we have to verify that
for any Borel set $B\in \mathscr F_T$, the intersection $B\cap S_n$ is a member of Borel  field $\mathscr F_{S_n}$ generated by open sets of the metric space $(S_n,\varrho_n)$. By Lemma \ref{lem1}, for any open subset $B$ of $T$,
$B\cap S_n$ is an open subset of $S_n$. Therefore, $\{ B\cap S_n; B\in \mathscr F_T\}$ is a sub-$\sigma$ field of $\mathscr F_{S_n}$ and thus the definition of $\tilde \mu_n$ is justified.

To apply  Dudley's theorem, we also need to verify that $\tilde\mu_n\stackrel{w}\to \tilde\mu_0$. To this end, we only need to verify that
\begin{equation}
\liminf_{n\to\infty}\tilde\mu_n(B)\ge\tilde\mu_0(B),\ \ \forall B\in\mathscr O_T,
\label{eqweak}
\end{equation}
where $\mathscr O_T$ is the collection of all open subsets of $(T,\de)$.
By assumption  $\mu_n\varphi_n^{-1}\stackrel{w}\to \mu_0$ and Theorem 2.1 of \cite{Billingsley99C}, we have
\begin{eqnarray}
\liminf_{n\to\infty}\tilde\mu_n(B)&=&\liminf_{n\to\infty}\mu_n(B\cap S_n)\nonumber\\
&\ge& \liminf_{n\to\infty}\mu_n(\varphi_n^{-1}(B_0))\ge \mu_0(B_0),
\label{eqweak0}
\end{eqnarray}
where the second inequality follows by applying Lemma \ref{lem2} and third by assumption.

Therefore, by  Dudley's theorem, there is a sequence of random elements $\tilde X_n$ with distributions $\tilde\mu_n$ defined on a common  probability space $(\Omega, \mathscr F,P)$ and  such that
$\tilde X_n\stackrel{\de}\to \tilde X_0, a.s.$ By Lemma \ref{lem3}, we conclude that $\varphi_n(X_n)\stackrel{\varrho_0}\to X_0$ a.s..

What remains to show is the case where the mapping $\varphi_n$ is measurable but not continuous. In this case, applying Luzin's theorem, for each $n>0$, we can find an continuous mapping $\tilde \varphi_n$ such that
\begin{equation}\label{eqluzin}
\mu_n\big(y\in S_n;\, \varphi_n(y)\ne \tilde \varphi_n(y)\big)<2^{-n}.
\end{equation}
By what we have proved for the case of continuous mappings, there exist $S_n$-valued random elements $X_n$ defined on a common probability space $(\Omega,\mathscr F,P)$ and the  distribution $X_n$ is $\mu_n$ and satisfy $\varrho_0(\tilde \varphi_n(X_n),X_0)\to 0, a.s.$ By (\ref{eqluzin}), we have
$$
\sum_{n=1}^\infty P\big( \varphi_n(X_n)\ne \tilde \varphi_n(X_n)\big)\le 1.
$$
By Borel Cantelli lemma, we conclude that the sequences $\{\varphi_n(X_n)\}$ and $\{\tilde \varphi_n(X_n)\}$ converge simultaneously with probability 1.
   The  proof  of  the  theorem  is
complete. $\hfill{} \Box$

\begin{remark}
  In  general  the  space  $T$  is  not  complete  under  $\de$, even  if  all  $S_n$   are  under  $\varrho_n$.   To  see  this, consider   the   case   that   all   $x_n(n>0)$    lie   in   $S_m$    for   one   fixed   $m$.    Then   $\{x_n\}$   is   $\de$-Cauchy   iff $\{(x_n,\varphi_m(x_n))\}_{n>0}$ is  $\varrho_m\times\varrho_0$-Cauchy.   If the  latter  holds,  then  $\{(x_n,\varphi_m(x_n))\}$  converges  in  $S_m\times S_0$,
but  not  necessarily  in  graph $\varphi_m$,  unless  the  latter  is  closed.   This  combined  with  the  observation  that
$\de$-Cauchy sequences $\{x_n\}_{n >0}$  with $x_n \in  S_n$  converge  if  $S_0$ is  $\varrho_0$-complete  leads  us  to  the  following  result.
\end{remark}

\begin{theorem}\lb{th3}
   Let  $S_n$   be  separable  and   $\varrho_n$-complete  for  each  $n$.   Then  $T$   is  $\de$-complete  iff  graph $\varphi_n$  is
closed in  $S_n\times S_0$ for  each  $n$.
\end{theorem}

It is  well-known  that  graph $\varphi_n$  is  closed  if $\varphi_n$ is  continuous,  and  that $\varphi_n$ is  continuous  if  graph $\varphi_n$  is
closed  and  $S_0$ is  compact. If a set $G$ is
the intersection of at most countably many open sets, then $G$  is called a $G_\de$. Using  the  fact  that  a  subset  of a  Polish  space  is  Polish  iff  it  is  $G_\de$   (see  Theorem 8.3 in Chapter XIV of \cite{Dugundji66T})),  we  arrive  at  the  following  variation  on Theorem \ref{th3}.

\begin{theorem}
   Let $S_n$  be Polish for  each  $n$.   Then  $T$  is  Polish  iff  graph $\varphi_n$ is  $G_\de$ in  $S_n\times S_0$ for  each  $n$.
\end{theorem}

\section{Some Lemmas}
In this section, we present some basic lemmas which are used to prove Theorem \ref{th1}.
\begin{lemma}
\label{lem1}
If $x\in S_n$ is an inner point of $B$, a subset of the metric space $(T,\de)$ and $\varphi_n$ is continuous with the metric $\varrho_n$, then $x$ is an inner point of
$B\cap S_n$, a subset of the metric space $(S_n,\varrho_n)$.
\end{lemma}

\noindent {\bf Proof.} Since $x$ is an inner point of $B$, there is an $r>0$ such that $B_{\de}(x,r)\subset B$. By continuity of $\varphi_n$, there is a positive constant
$\eta$ such that for any $y\in B_{\varrho_n}(x,\eta)$ we have $\varrho_0(\varphi_n(x),\varphi_n(y))\le r/2$. Thus, for any $y\in B_{\varrho_n}(x,\eta)\cap B_{\varrho_n}(x,r/2)$,
we have
$$
\de(x,y)=\varrho_0(\varphi_n(x),\varphi_n(y))+\ep_n\wedge\varrho_n(x,y)<r,
$$
Thus $y\in B_{\de}(x,r)$. Noting that $B_{\varrho_n}(x,\eta)\cap B_{\varrho_n}(x,r/2)$ is an open subset of $B_{\de}(x,r)\cap S_n$, then the proof of the lemma is complete.

\begin{lemma}
\label{lem2}
If $B$ is an open set of $(T,\de)$ and $B_0=B\cap S_0$, then $\varphi_n^{-1}B_0\subset B\cap S_n$ for all large $n$.
\end{lemma}

\noindent {\bf Proof.} Suppose $x_0\in B_0$. Since $x_0$ is an inner point of $B$, there is an open ball
$B_{\de}(x_0,r)\subset B$. Now, assume that $n$ is so large that $\ep_n<r$. If $x\in S_n$ is such that $\varphi_n(x)=x_0$,
then
$
\de(x_0,x)=\ep_n<r,
$
which implies  $x\in B_{\de}(x_0,r)$ and consequently, $\varphi_n^{-1}(x_0)\subset B\cap S_n$. The proof is complete.

\begin{lemma}
\label{lem3}
With the metric defined in (\ref{eq1}),  if $x_n\stackrel{\de}\to x_0$ and $x_0\in S_k$ ($k>0$), then for almost all $n$, $x_n\in S_k$ and
 $\varrho_k(x_n,x_0)\to 0$.

 If $x_n\stackrel{\de}\to x_0$ and $x_0\in S_0$, then $s(x_n)\to\infty$ and $\varphi_{s(x_n)}(x_n)\stackrel{\varrho_0}\to x_0$.
\end{lemma}

\noindent {\bf Proof.} If $x_0\in S_k$ and there are infinitely many $n$ such that $s(x_n)\ne k$, then there are infinitely many $n$ such that
$\de(x_n,x_0)\ge \ep_k>0$ which violates to the assumption that $\de(x_n,x_0)\to 0$. Therefore, for almost all $n$, $s(x_n)=k$. Thus, $\varrho_k(x_n,x_0)\to 0$
follows from the simple fact that $\de(x_n,x_0)\ge \ep_k\wedge \varrho_k(x_n,x_0)$.

If $x_0\in S_0$ and there are infinitely  many $n$ such that $s(x_n)\le N$, then there are infinitely many $n$ such that
$\de(x_n,x_0)\ge \min_{k\le N}\ep_k>0$. Therefore, we have $s(x_n)\to \infty$. Thus, for all large $n$, we have
$$
\de(x_n,x_0)\ge\varrho_0(\varphi_{s(x_n)}(x_n),x_0)\to 0.
$$
The proof is complete.

\section{Applications}
To begin with, there is one of the simplest example which can be proved by our  Theorem \ref{th1}, but not by the
theorem of Skorohod-Dudley in its original form. It is Theorem 3.1 in \cite{Billingsley99C}, restricted to
separable metric spaces.

\begin{example}{\expandafter{\rm
  If $S$ is a separable metric space with metric $\varrho$, $(X_n, Y_n)$ are $S^2$-valued random variables for
$n = 1,2, \dots$ and $X$ is an $S$-valued random variable such that $X_n \toD X$ in $S$ and $\varrho(X_n , Y_n) \toD 0$ in $\mathbb{R}$,
then $Y_n \toD X$ in $S$.}}
\end{example}
\begin{proof}
  By Theorem 3.9 in \cite{Billingsley99C} we have $(X_n ,\varrho (X_n,Y_n)) \toD (X ,0)$ in $S \times \mathbb{R}$. Applying  Theorem \ref{th1}
with $S_0 = S \times \mathbb{R}$, $S_n = S^2$, $X_n$ replaced by $(X_n, Y_n)$ and $(\varphi_n(x,y) = (x ,
\varrho(x,y))$,  we obtain this conclusion.
\end{proof}
Next we will give some applications of the  strong representation theorem in
random matrix theory.
 Our Theorem \ref{th1} has a wide range of applications in random matrices,  especially in their spectral properties.  Next we give several examples. Before that we introduce some basic definitions.

\begin{definition}[ESD]
For any $n\times n$ matrix $\A$ with real eigenvalues, we define the  empirical spectral distribution (ESD) of $\A$ by
\begin{align*}
    F^{\A}(x)=\frac{1}{n}\sum_{i=1}^nI(\la_i^{\A}\le x),
\end{align*}
where $\la_i^{\A}$ is the $i$-th smallest eigenvalue of $\A$ and $I(B)$  is the indicator function of an event $B$.
\end{definition}
\begin{definition}[LSD]
  Let $\{\A_n\}$ be a sequence of random matrices with ESD $F^{\A_n}$. If $F^{\A_n}$ has a limit distribution $F$, then $F$  is called the limiting empirical distribution (LSD) of the sequence $\{\A_n\}$.
\end{definition}
\begin{definition}[Wigner matrix]
 Suppose $\W_n=n^{-1/2}\(w_{ij}\)_{i,j=1}^n$ is a  Hermitian matrix whose entries are all zero-mean
random variables. Then $\W_n$ is said to be a Wigner matrix if the following conditions are satisfied:
\begin{itemize}
  \item $\{w_{ij};1\leq i\leq j\leq n\}$are independent random variables;
  \item $E|w_{ij}|^2=1$, for all $1\leq i<j\leq n$.
\item For any $\eta>0$, as  $n\to\infty$,
\begin{align*}
    (\eta\sqrt{n})^{-2}\sum_{i,j}E(|w_{ij}|^2I(|w_{ij}|\geq \eta\sqrt n))\to0.
\end{align*}
\end{itemize}
\end{definition}
\begin{definition}[Stieltjes transform]
  For any function of bounded variation $H$ on the real line, its Stieltjes transform is defined by
$$s_H(z)=\int\frac{1}{\lambda-z}dH(\lambda),~~z\in\mathbb{C}^{+}\equiv\{z\in\mathbb{C}^{+}:\Im z>0\}.$$
\end{definition}
Then we have the following examples:

\begin{example}{\expandafter{\rm
Let $\X_1,\cdots,\X_n$ be an i.i.d. sample from a $d$-dimensional normal distribution with mean vector $\bbmu$ and covariance $\bbgS=\bbV\diag[\gl_1\bbI_{d_1},\cdots,\gl_k\bbI_{d_k}]\bbV'$,
where $\gl_1>\gl_2>\cdots>\gl_k\ge 0$ are distinct eigenvalues with multiplicities $d_1,\cdots, d_k$ ($d_1+\cdots+d_k=d$) of the population covariance matrix $\bbgS$,  and $\bbV$ is an orthogonal matrix of orthonormal eigenvectors of $\bbgS$.
Write $\S_n=\frac1{n-1}\sum_{i=1}^n(\X_i-\bar \X)(\X_i-\bar \X)'$ be the sample covariance matrix. By the law of large numbers, we have $\S\to\bbgS, a.s.$ and hence
$$
l_t\to\gl_j, \mbox{ if } d_1+\cdots+d_{j-1}<t\le d_1+\cdots+d_j,
$$
where $l_1\ge l_2\ge \cdots\ge l_d$  are the ordered eigenvalues of $\S_n$.

Now, we investigate the limiting distribution of
$$\{\sqrt{n}(l_t-\gl_j), d_1+\cdots+d_{j-1}<t\le d_1+\cdots+d_j, j=1,\cdots,k\}.$$
To begin with, we consider the limiting distribution of $\M_n=\sqrt{n}(\S_n-\bbgS)$. 
By classical CLT, it is easy to see that $\M_n$ tends to a $d\times d$ symmetric random matrix $\M=(m_{ij})$ in distribution, where $m_{ij}, i\le j$
are jointly normally distributed with means 0 and covariances
$$
{\rm Cov}(m_{ij}, m_{ts})={\rm Cov}(X_{i1}X_{j1}, X_{s1}X_{t1}).
$$

 Define a measurable mapping $\varphi_n$ from $\mathbb R^{d\times n}$ to $\mathbb R^{\frac12 d(d+1)}$ such that
$$
\varphi_n(\X)=\M_n=\sqrt{n}(\S_n-\bbgS).
$$
Applying Theorem \ref{th1}, we may redefine $\widetilde\X=(\widetilde\X_1,\cdots,\widetilde\X_n)$ on $\mathbb R^{d\times n}$ and $\widetilde\M$ on $\mathbb R^{\frac12d(d+1)}$ on a suitable probability space satisfying
$\widetilde\M_n\to \widetilde\M$ a.s.. Blocking the matrices $\bbV=(\bbV_1,\cdots,\bbV_k)$, $\bbV'\widetilde\M_n\bbV=(\widehat\M_{n,ij}),$ and $\bbV'\widetilde\M\bbV=(\widehat\M_{ij})$ where $\bbV_j$ consists of the $d_j$ eigenvectors of $\gl_j$,
$\widehat\M_{n,ij}=\bbV_i'\widetilde\M_n\bbV_j$ and $\widehat\M_{ij}=\bbV_i'\widetilde\M\bbV_j$.

Denote the spectral decomposition of $\widetilde\S_n=\U_n\diag[\widetilde l_1,\cdots,\widetilde l_d]\U_n'$ and split the matrices as blocks $\U_n=(\U_{n,1},\cdots,\U_{n,k})$ and $\diag[\widetilde l_1,\cdots,\widetilde l_d]=\diag[\D_{n,1},\cdots,\D_{n,k}]$ accordingly.
Then, $\widetilde\M_n\to \widetilde\M, a.s.$ is equivalent to
\begin{gather}\label{eqb1}\sqrt{n}\Big(\bbV'\U_n\diag[\D_1,\cdots,\D_k]-\diag[\gl_1\I_{d_1},\cdots,\gl_k\I_{d_k}]\bbV'\U_n\Big)\\
-(\widehat\M_{ij})\bbV'\U_n\to 0,\quad a.s..\nonumber
\end{gather}
The $(i,j)$-block with $i\ne j$ of (\ref{eqb1}) is
$$
\sqrt{n}\bbV_i'\U_{n,j}(\D_j-\gl_i\I_{d_j})-\sum_{t=1}^k\widehat\M_{it}\bbV_t'\U_{n,j}\to 0,~~a.s.
$$
which together with the fact that $\D_j-\gl_i\I_{d_j}\to(\gl_j-\gl_i)\I_{d_j}$ and $\gl_i\ne\gl_j$ implies  $\bbV_i'\U_{n,j}=O(1/\sqrt{n})$.
Consequently, we obtain that
$$
\bbV_i'\U_{n,i}\U_{n,i}'\bbV_i=\I_{d_i}-\sum_{t\ne i}\bbV_i'\U_{n,t}\U_{n,t}'\bbV_i=\I_{d_i}+O(1/n).
$$
This proves that $\bbV_i'\U_{n,i}$ is asymptotically orthogonal. What is more, the $(i,i)$ block of (\ref{eqb1}) is
\begin{gather*}
\bbV_i'\U_{n,i}\left(\sqrt{n}(\D_i-\gl_i\I_{d_i})\right)-\sum_{t=1}^k\widehat\M_{it}\bbV_t'\U_{n,i}\\
=\bbV_i'\U_{n,i}\left(\sqrt{n}(\D_i-\gl_i\I_{d_i})\right)-\widehat\M_{ii}\bbV_i'\U_{n,i}+o(1)\to 0, ~~a.s..
\end{gather*}
Therefore, $\sqrt{n}(\D_i-\gl_i\I_{d_i})$ tends to a diagonal matrix of ordered eigenvalues of the matrix $\widehat\M_{ii}$ and $\bbV_i'\U_{n,i}$  tends to the matrix of orthonormal eigenvectors of $\widehat\M_{ii}$ if we suitably
select the signs of the eigenvectors.

Checking the covariances of the entries of $\widetilde\M$, we have
the variances of the diagonal entries of $\widehat\M_{ii}$ is $2\gl_i^2$ and that of off-diagonal elements are $\gl_i^2$. Also, the covariances of the entries $\widehat\M_{ii}$ and that of $\widehat\M_{jj}$ $(i\ne j)$ are 0.

Therefore, we conclude that the random vector $\{\sqrt{n}(\D_j-\gl_j\I_{d_j}),\ j=1,\cdots, k\}$ tends to $k$ independent sub-vectors and its $j$-th sub-vector consists of the ordered eigenvalues of a Wigner matrix
whose diagonal entries are $N(0,2\gl_j^2)$ and off-diagonal entries are $N(0,\gl_j^2)$.}}
 \end{example}
\begin{remark}
The random vectors may not be necessarily assumed normal. Under certain moment assumptions, the result remains true.
\end{remark}

\begin{remark}
Anderson in \cite{anderson51} considered the limiting distributions of the relative eigenvalues of two independent Wishart matrices without using
the strong representation theorem. As a consequence, he has to argue the continuity of the inverse transformation of spectral decomposition. In fact,
the inverse transformation is not completely continuous, it has a lot of exception points and it is easy to argue that the exception points form
a set of probability zero. Using the strong representation theorem. We do not need to worry about the probability of exception points.
\end{remark}

\begin{remark}
When dimension $d$ is fixed, the result can also be proved by using the original version of Skorokhod strong representation theorem. In this case the metric space can be chosen as $\mathbb R^{\frac12d(d+1)}$ and the random elements  are $\widetilde\M_n=\sqrt{n}(\widetilde\S_n-\bbgS)$, where $(n-1)\S_n$ is a Wishart random matrix. Then the derivation will be the same as above. However, when normality is not assumed, the structure of
sample covariance matrix of $\widetilde\S_n$ will be lost.

Furthermore, if the dimension of the population increases as the sample size increases, the original version of Skorokhod strong representation
theorem is not applicable. See the next example.
\end{remark}

\begin{example}{\expandafter{\rm Silverstein  in \cite{jack95} proved the following result. Let $\T^{1/2}_n$ be the Hermitian non-negative square root of a $p\times p$ $\T_n$, and let $\B_n =\frac1n \T^{1/2}_n \X_n\X^*_n\T^{1/2}_n$, where the  $\T_n$ is independent of $\X_n$ and its ESD almost surely  tends to a proper cumulative distribution function (c.d.f.), $\X_n$ is a $p\times n$ matrix whose entries are i.i.d. random variables with mean 0 and variance 1, and $p/n\to y>0$. Then, almost surely, $F^{\B_n}$ converges in distribution, as $n \to\infty$, to a (nonrandom) c.d.f. $F$, whose Stieltjes transform
$m(z)$ ($z \in \mathbb C^+$) satisfies
$$
m =\int\frac{1}{\tau(1 -y- yzm) - z}dH(\tau):
$$
in the sense that, for each $z \in \mathbb C^+$, $m = m(z)$ is the unique
solution to te equation above in $D_y=\{ m\in\mathbb C : −\frac{(1-y)}z + ym \in \mathbb C^+\}$.

Now, we want to show that if the ESD of $\T_n$ tends to a proper c.d.f. $H$ in probability, the result remains true provided to weaken the strong convergence of $F^{\B_n}$ to convergence in probability. Applying Theorem \ref{th1} with $S_n=\mathbb R^{pn+\frac12p(p+1)}$ with random elements
$\{(\X_n,\T_n)\}$, $\varphi_n(\X_n,\T_n)=F^{\T_n}$, $S_0$ as the collection of c.d.f. the limiting element $H$, then we can construct a probability space
$(\Omega,\mathscr F, P)$ on which we have $(\widetilde \X_n,\widetilde \T_n)$ with identical distributions as $(\X_n, \T_n)$ and satisfies
$F^{\widetilde \T_n}\to H, a.s.$. Then, applying the results of Silverstein \cite{jack95}, we obtain
$$
F^{\widetilde \B_n}\to F, a.s.
$$
where $\widetilde \B_n=\frac1n \widetilde \T_n^{1/2}\widetilde \X_n\widetilde \X_n^*\widetilde \T_n^{1/2}$. Because $F^{\B_n}\stackrel{\mathcal{D}}=F^{\widetilde \B_n}$, we conclude that
$$
F^{ \B_n}\stackrel{{p}}{\to} F.
$$

Note that $\T_n$ and $H$ do not tack values in a common metric space, the original version of Skorokhod theorem is not applicable.}}
\end{example}

Similarly, due to Theorem 1.1 in \cite{BaiZ10L} we obtain the following result.

\begin{example}{\expandafter{\rm
 For each $n = 1, 2, \dots$, let $\W_n$ be a Wigner matrix as defined above and let $\T_n$ be a
Hermitian nonnegative definite matrix with $(\T^{1/2})^2=\T$. Suppose that, as $n\to\infty$, the empirical spectral distribution
of $\T_n$ converges weakly to a non-random probability distribution $H$ in probability. Let $\B_n = n^{-1/2} \T^{1/2}_n \W_n \T^{1/2}_n$. Then, as $n\to\infty$, the ESD of  $\B_n$ converges weakly to a non-random probability
distribution $F$ in probability, whose Stieltjes transform $s(z)$ uniquely solves the following equation system
\begin{align*}
    \left\{
       \begin{array}{ll}
         s(z)=-z^{-1}-z^{-1}(g(z))^2,  \\
         g(z)=\int\frac{t}{-z-tg(z)}dH(t),
       \end{array}
     \right.
\end{align*}
for any $z\in \mathbb{C}^+=\{z\in \mathbb{C}: \Im z>0\}$, where $g(z) \in \mathbb{C}$ with $\Im
 g(z) \geq 0$.}}
\end{example}

\begin{example}\lb{th6}{\expandafter{\rm
  Let $\mathcal{C}$ be a connected open set of complex plane $\mathbb{C}$ and $\{Y_n(z),z\in\mathcal{C}\}$ be a  two-dimensional stochastic process which is defined on a separable metric space $S_n,n=0,1,2,\dots$. Suppose that $Y_n(z)$ is analytic and bounded by a constant for every $n$ and $z\in\mathcal{C}$. If as $n\to\infty$, $Y_n(z)$ converges weakly to $Y_0(z)$ for
each $z$ in a subset of $\mathcal{C}$. Then we have $Y'_n(z)$ converges weakly to $Y'_0(z)$  for all $ z\in\mathcal{C}$, where $'$ denote the derivative of the function $Y_n$ at $z$, $n=1,2,\dots,\infty$.}}
\end{example}
\begin{proof}
  Applying Theorem \ref{th1}, we get that there is one probability space on
which we can define a two-dimensional stochastic process $\{\hat Y_n(z),z\in\mathcal{C}\}$, such that, for
each $n=1,2,\dots$, the  distribution of $\{\hat Y_n(z),z\in\mathcal{C}\}$ is identical to that of $\{Y_n(z),z\in\mathcal{C}\}$ and
 $\{\hat Y_n(z),z\in\mathcal{C}\}$ converges to $\{Y_0(z),z\in\mathcal{C}\}$ almost surely for
each $z$ in the subset of $\mathcal{C}$. Then using Vitali's convergence theorem (see Lemma 2.3 in \cite{BaiS04C}),
 we  obtain  that $\hat Y'_n(z)$
 converge  almost  surely  to $\hat Y'_0(z)$
 for  all $z\in\mathcal{C}$, which implies $Y'_n(z)$ converges weakly to $Y'_0(z)$. The proof of this example  is complete.
\end{proof}

Combining Example \ref{th6} and Theorem 2.1 in \cite{BaiY05C}, we can get  the following conclusion:
\begin{example} Let $F^{sc}$ be the LSD of Wigner matrices  $\{\W_n\}$.
Suppose that:
\begin{itemize}
  \item[(\romannumeral1)] For all $i$, $E|w_{ii}|^2=\sigma>0$, and if $\W_{n}$ is complex, $Ew_{ij}^2=0$ for all $i<j$.
  \item[(\romannumeral2)]   $E|w_{ij}|^4=M\leq \infty$, $i\neq j$;
  \item[(\romannumeral3)] For any $\eta>0$, as  $n\to\infty$,
\begin{align*}
    (\eta\sqrt{n})^{-4}\sum_{i,j}E(|w_{ij}|^4I(|w_{ij}|\geq \eta\sqrt n))\to0.
\end{align*}
\end{itemize}
Then we get that  the process $\{n[s_{F^{\W_n}}'(z)-s_{F^{sc}}'(z)]; z\in \mathcal{C}_{sc}\}$ converges weakly to a
Gaussian process $\{\xi(z);\mathcal{C}_{sc}\}$ with the mean and covariance functions:
\begin{align*}
    E\xi(z)=a'(z)\qquad\mbox{and}\qquad Cov(\xi(z_1),\xi(z_2))=\frac{\partial^2 b(z_1,z_2)}{\partial z_1\partial z_2},
\end{align*}
where
\begin{gather*}
\mathcal{C}_{sc}=\{z=u+iv|u\not\in[-2,2],|v|\geq v_0>0\},\\
a(z)=[1+s_{F^{sc}}'(z)]s_{F^{sc}}^3(z)[\sigma^2-1+(\kappa-1)s_{F^{sc}}'(z)+\b s_{F^{sc}}^2(z)],\\
b(z_1,z_2)=s_{F^{sc}}'(z_1)s_{F^{sc}}'(z_2)[\sigma^2-\kappa+2\b s_{F^{sc}}(z_1)s_{F^{sc}}(z_2)+{\kappa}{(1-s_{F^{sc}}(z_1)s_{F^{sc}}(z_2))^{-2}}],\\
    s_{F^{sc}}(z)=-\frac{1}{2}\(z-\sqrt{z^2-4}\),\\
\b=E(|x_{12}|^2-1)^2-\kappa\quad\mbox{and}\quad\kappa=\left\{
                                            \begin{array}{ll}
                                              1, & \hbox{$\W_n$ is  complex;} \\
                                              2, & \hbox{$\W_n$ is real.}
                                            \end{array}
                                          \right.
\end{gather*}
\end{example}

%

\bibliographystyle{abbrv}

\begin{thebibliography}{10}

\bibitem{anderson51}

Anderson, T. W. (1951). The asymptotic distribution of certain characteristic roots and vectors,
{\it Proc. 2nd Berkeley Symp. Math. Statist. Prob.}, pp 103--130.

\bibitem{Bailiang}
Z.~D. Bai and W. Q. Liang.
\newblock {Strong representation of week convergence}.
\newblock {\em Technical rept.  No. 85-29,} Pittsburgh University, Center for Multivariate Analsis, 1985.

\bibitem{Bailv}
Z.~D. Bai, W. Q. Liang, and W. Vervaat.
\newblock {Strong representation of week convergence}.
\newblock {\em Technical rept. Sep 86-Sep 87,} North Carolina Univ At Chapel Hill, Dept Of Statistics, 1987.

\bibitem{BaiS04C}
Z.~D. Bai and J.~W. Silverstein.
\newblock {CLT for linear spectral statistics of large-dimensional sample
  covariance matrices}.
\newblock {\em The Annals of Probability}, 32(1):553--605, 2004.

\bibitem{BaiY05C}
Z.~D. Bai and J.~F. Yao.
\newblock {On the convergence of the spectral empirical process of Wigner
  matrices}.
\newblock {\em Bernoulli}, 11(6):1059--1092, Dec. 2005.

\bibitem{BaiZ10L}
Z.~D. Bai and L.~X. Zhang.
\newblock {The limiting spectral distribution of the product of the Wigner
  matrix and a nonnegative definite matrix}.
\newblock {\em Journal of Multivariate Analysis}, 101(9):1927--1949, Oct. 2010.

\bibitem{BertiP10S}
P.~Berti, L.~Pratelli, and P.~Rigo.
\newblock {A Skorohod representation theorem for uniform distance}.
\newblock {\em Probability Theory and Related Fields}, 150(1-2):321--335, Mar.
  2010.

\bibitem{Billingsley99C}
P.~Billingsley.
\newblock Convergence of probability measures, 2nd edition.
\newblock Wiley Series in Probability and Statistics, 1999.

\bibitem{BlackwellD83E}
D.~Blackwell and L.~Dubins.
\newblock {An extension of Skorohod's almost sure representation theorem}.
\newblock {\em Proc. Amer. Math. Soc}, 89(4):691--692, 1983.

\bibitem{Dudley68D}
R.~M. Dudley.
\newblock {Distances of probability measures and random variables}.
\newblock {\em The Annals of Mathematical Statistics}, 39(5):1563----1572,
  1968.


\bibitem{Dugundji66T}
J. Dugundji.
\newblock {Topology}.
\newblock { Allyn and Bacon},
  1966.

\bibitem{Sethuraman02S}
J.~Sethuraman.
\newblock {Some extensions of the Skorohod representation theorem}.
\newblock {\em Sankhya: The Indian Journal of Statistics, Series A},
  64(3):884--893, 2002.

\bibitem{jack95}

Silverstein, J.W. (1995). Strong convergence of the empirical distribution of eigenvalues
of large dimensional random matrices {\it J. Multivariate Anal.} {\bf 55} 331--339.


\bibitem{Skorokhod56L}
A.~V. Skorokhod.
\newblock Limit theorems for stochastic processes.
\newblock {\em Theory of Probability $\&$ Its Applications}, 1(3):261--290,
  1956.

\bibitem{Whitt02S}
W.~Whitt.
\newblock {\em Stochastic-process limits: an introduction to stochastic-process
  limits and their application to queues}.
\newblock Springer, 2002.

\end{thebibliography}

\end{document}